\documentclass[10pt]{article}
\usepackage{amsmath,amssymb,amsbsy,amsfonts,latexsym,amsthm,amscd,amsxtra,
amssymb,appendix}
\usepackage[all]{xy}

\newtheorem{Teo}{Theorem}[section]
\newtheorem{Def}[Teo]{Definition}
\newtheorem{Cor}[Teo]{Corollary}
\newtheorem{Lem}[Teo]{Lemma}
\newtheorem*{acknow}{Acknowledgements}

\newtheorem{Prop}[Teo]{Proposition}

\newtheorem{Rema}[Teo]{Remark}

\newtheorem*{Pf}{Proof}
\newenvironment{Proof}{\begin{Pf} \begin{upshape}} {\end{upshape} \qed\end{Pf}}

\newcommand\beq[1]{ \begin{equation} \label{#1}}
\newcommand{\eeq}{ \end{equation} }
\newcommand{\beqno}{ \begin{equation*}}
\newcommand{\eeqno}{ \end{equation*}}
\newcommand\beqa[1]{ \begin{eqnarray} \label{#1}}
\newcommand{\eeqa}{ \end{eqnarray} }
\newcommand{\beqano}{ \begin{eqnarray*} }
\newcommand{\eeqano}{ \end{eqnarray*} }

\def\be{\begin{equation}}
\def\ee{\end{equation}}
\def\bea{\begin{eqnarray}}\def\eea{\end{eqnarray}}
\def\nn{\nonumber}

\newcommand{\T}{ {\mathbb T}   }
\newcommand{\N}{ {\mathbb N}   }
\newcommand{\R}{ {\mathbb R}   }
\newcommand{\Z}{ {\mathbb Z}   }

\def\TTT{\T}%{\hbox{\msytw T}}
\def\RRR{\R}%{\hbox{\msytw R}}
\def\ZZZ{\Z}%{\hbox{\msytw Z}}

\renewcommand \a {\alpha}
\newcommand \e {\varepsilon }

\renewcommand \b  {\beta}
\renewcommand \r {\rho}
\renewcommand \d {\delta}

\newcommand \m {\mu}
\newcommand \n {\nu}

\newcommand \f {\varphi}

\newcommand \g {\gamma}

\renewcommand \t {\tau}

\renewcommand \l {\lambda}
\renewcommand \L {\Lambda}
\newcommand \ph {\phi}

\newcommand \p {\pi}

\newcommand \cA {{\mathcal A}}

\newcommand \cS {{\mathcal S}}

\newcommand \cL {{\mathcal L}}

\newcommand \dpr {\partial}

\newcommand \rH {{\rm H}}

\newcommand \calM {\mathfrak{M}}
\newcommand \cM {{\mathcal M}}

\def\TT{{\mathcal T}}
 
\def\LL{{\mathcal L}}  
\def\DD{{\mathcal D}}

\let\0=\noindent

\def\media#1{{\langle#1\rangle}}
\def\ie{\hbox{\it i.e.,\ }}

\let\dpr=\partial

\let\o=\omega
\let\==\equiv

\begin{document}

\title{Uniqueness of Invariant Lagrangian Graphs in a Homology or 
a Cohomology Class.}
\author{\sc{Albert Fathi, Alessandro Giuliani and Alfonso Sorrentino}}

\date{}
\maketitle

%%%%%%%%%%%%%%%%%%%%%%%%%%%%%%%%%%%%%%%%%%%%%%%%
%%%%%%%%%%%%%%%%%% ABSTRACT %%%%%%%%%%%%%%%%%%%%
\begin{center}
\parbox[c]{4.5 in}
{\small Given a smooth compact Riemannian manifold $M$ and a Hamiltonian $H$
on the cotangent space ${\rm T}^*M$, strictly convex and superlinear in the 
momentum variables, we prove uniqueness of certain ``ergodic'' invariant 
Lagrangian graphs within a given homology or cohomology class. In particular, 
in the context of quasi-integrable Hamiltonian systems, our result implies 
global uniqueness of Lagrangian KAM tori with rotation vector $\r$. This result
extends generically to the $C^0$-closure of KAM tori.}
\end{center}

%%%%%%%%%%%%%%%%%%%%%%%%%%%%%%%%%%%%%%%%%%%%%%%%%%%%%%%%%%%%%%%%%%%%%%%%%%
%%%%%%%%%%%%%%%%%%%%%%%% SECTION 1%%%%%%%%%%%%%%%%%%%%%%%%%%%%%%%%%%%%%%%%
%%%%%%%%%%%%%%%%%%%%%%%%%%%%%%%%%%%%%%%%%%%%%%%%%%%%%%%%%%%%%%%%%%%%%%%%%%

\section{Introduction}\label{intro}

%%%%%%%%%%%%%%%%%%%%s+\t}%%%%%%%%%%%%%%%%%%%%%%%%%%%%%%%%%%%%%%%%%%%%%%%%%%%%%%
%%%%%%%%%%%%%%%%%%%%%%%%%%%%%%%%%%%%%%%%%%%%%%%%%%%%%%%%%%%%%%%%%%%%%%%%%%

A particularly interesting and fruitful approach to the study of local and 
global properties of dynamical systems is concerned with the study of invariant
submanifolds, rather than single orbits, paying particular 
attention to their existence, their ``fate'' and their geometric properties.
In the context of quasi-integrable Hamiltonian systems, one of the most 
celebrated breakthroughs in this kind of approach was KAM theory, which 
provided a method to construct invariant submanifolds diffeomorphic 
to tori, on which the dynamics is conjugated to a quasi-periodic motion
with rotation vector $\r$, sometimes referred to as {\it KAM tori}. 
KAM theory finally settled the old question about {\it existence} of such 
invariant submanifolds 
in ``generic'' quasi-integrable Hamiltonian systems, dating back at least
to Poincar\'e, and opened the way to a new understanding of the nature of 
Hamiltonian systems, of their stability and of the onset of chaos in classical 
mechanics. However, the natural question about the {\it uniqueness} of 
these invariant submanifolds for a fixed rotation vector $\r$
remained open for many more years and, quite 
surprisingly, even nowadays, for many respects it is still unanswered. 
A possible reason for this is that the analytic methods, which the 
KAM algorithm is based on, are not well suited for studying global
questions, while, on the other hand, the natural variational methods to 
approach this problem have been developed only much more recently and 
they are still not so widely well-known.

In this paper we address the above problem and prove global uniqueness of 
certain ``ergodic'' invariant Lagrangian graphs 
within a given homology or cohomology class, for a large class 
of convex Hamiltonians, known as Tonelli Hamiltonians.
Our work will be based on the variational 
approach provided by the so-called Aubry-Mather theory, as well as its 
functional counterpart called weak KAM theory.

The paper is organized as follows. In Section \ref{sec1} we  
define the geometric objects we shall look at (the invariant Lagrangian 
graphs), we introduce some concepts (homology class of an 
invariant measure and Schwartzman ergodicity), which will turn out to be useful
for illustrating their dynamical properties, and we state our main 
uniqueness results.
In Section \ref{sec3} 
we prove our main results. In Section \ref{sec4} we discuss %in detail
the implications of our results for KAM theory and compare
them with some previous local uniqueness theorems for KAM 
tori. 
In \ref{app1} we shall discuss some details
concerning the definition of Schwartzman ergodicity, give 
some examples and describe some properties of Schwartzman ergodic flows.

\vspace{.3truecm}

{\bf Acknowledgements.} We would like to thank Giovanni Gallavotti for 
having drawn our attention to this problem and for useful discussions. 
AF and AS are grateful to John Mather for many fruitful conversations.
%%%%%%%%%%%%%%%%%%%%%%%%%%%%%%%%%%%%%%%%%%%%%%%%%%%%%%%%%%%%%%%%%%%%%%%%%%
%%%%%%%%%%%%%%%%%%%%%%%% SECTION 1%%%%%%%%%%%%%%%%%%%%%%%%%%%%%%%%%%%%%%%%
%%%%%%%%%%%%%%%%%%%%%%%%%%%%%%%%%%%%%%%%%%%%%%%%%%%%%%%%%%%%%%%%%%%%%%%%%%

\section {Setting and Main results}\label{sec1}

%%%%%%%%%%%%%%%%%%%%%%%%%%%%%%%%%%%%%%%%%%%%%%%%%%%%%%%%%%%%%%%%%%%%%%%%%%
%%%%%%%%%%%%%%%%%%%%%%%%%%%%%%%%%%%%%%%%%%%%%%%%%%%%%%%%%%%%%%%%%%%%%%%%%%

Let $M$ be a compact and connected 
smooth manifold without boundary of dimension $n$.
Denote by ${\rm T}M$ its tangent bundle and ${\rm T}^*M$ the cotangent one. A
point of ${\rm T}M$ will be denoted by $(x,v)$, where $x\in M$ and $v\in
{\rm T}_xM$, and a point of ${\rm T}^*M$ by $(x,p)$, where 
$p\in {\rm T}_x^*M$ is a linear form on the vector space ${\rm T}_xM$. 
Let us fix a Riemannian metric $g$ on $M$ and let
$\|\cdot\|_x$ be the norm induced by $g$ on ${\rm T}_x M$; we shall use the 
same notation for the norm induced on the cotangent space
${\rm T}_x^*M$. This cotangent space ${\rm T}^*M$ can be canonically endowed 
with a symplectic structure, given by the exact $2$-form 
$\o=dx\wedge dp = -d(pdx)$, where $({\mathcal U},\,x)$ is a local coordinate 
chart for $M$ and $({\rm T}^*{\mathcal U},\,x, p)$ the associated 
cotangent coordinates. The $1$-form $\l = pdx$ is also called 
{\it tautological form} (or {\it Liouville form}) and is intrisically defined, 
\ie independently of the choice of local coordinates.
A distinguished role in the study of the geometry of a symplectic space is 
played by the so-called {\it Lagrangian submanifolds}.

\begin{Def}\rm
Let $\Lambda$ be an $n$-dimensional $C^1$ submanifold of 
$({\rm T}^*M,\omega)$. We 
say that $\Lambda$ is {\it Lagrangian} if for any $(x,p)\in \Lambda$, 
${\rm T}_{(x,p)}\Lambda$ is a Lagrangian subspace, \ie
$\omega \big|_{{\rm T}_{(x,p)}\Lambda} = 0$. 
\end{Def}

We shall mainly be concerned with Lagrangian graphs, that is Lagrangian 
manifolds $\Lambda\subset T^*M$ such that $\Lambda=\{(x,\eta(x))\,,\ x\in M\}$.
It is straightforward to check that the graph $\Lambda$ is Lagrangian if 
and only if $\eta$ is a closed $1$-form. 
The element $c=[\eta]\in \rH^1(M;\R)$ is called the {\it 
cohomology class}, or {\it Liouville class}, of $\Lambda$.
This motivates the following extension of the notion of 
Lagrangian graph to the Lipschitz case.

\begin{Def}\rm
A Lipschitz section $\L$ of ${\rm T}^*M$ is a {\it Lipschitz Lagrangian graph}
if it locally coincides with the graph of an exact differential.
\end{Def}

Observe that a Lipschitz Lagrangian graph $\Lambda$ is differentiable almost everywhere and at each differentiability point $(x,p)$ the tangent space ${\rm T}_{(x,p)}\Lambda$ is a Lagrangian subspace. Such Lipschitz graphs, although less regular, enjoy many properties of $C^1$-Lagrangian graphs. In the following,  when referring   to a Lagrangian graph without specifying its regularity, we shall assume that it is at least Lipschitz.\\

We shall consider the dynamics on ${\rm T}^*M$ generated by a Tonelli 
Hamiltonian.

\begin{Def}\rm A function $H:\,{\rm T}^*M\longrightarrow \R$ is called a 
{\it Tonelli (or optical) Hamiltonian} if:
\begin{itemize}
\item[i)]  the Hamiltonian $H$ is of class $C^k$, with $k\geq 2$;
\item[ii)]  the Hamiltonian $H$ is strictly convex in the fiber in the $C^2$ 
sense, \ie the second partial 
vertical derivative ${\dpr^2 H}/{\dpr p^2}(x,p)$ is positive definite, 
as a quadratic form, for any $(x,p)\in {\rm T}^*M$;
\item[iii)] the Hamiltonian $H$ is superlinear in each fiber, \ie
$$\lim_{\|p\|_x\rightarrow +\infty} \frac{H(x,p)}{\|p\|_x} = + \infty$$
(because of the compactness of $M$, this condition is independent of the choice
of the Riemannian metric).
\end{itemize}
\end{Def}

Given $H$, we can define the associated {\it Lagrangian},
as a function on the tangent bundle:
\bea L:\; {\rm T}M &\longrightarrow & \R \nn\\
(x,v) &\longmapsto & 
\sup_{p\in {\rm T}^*_xM} \{\langle p,\,v \rangle_x -H(x,p)\}\,\nn \eea 
where $\langle \,\cdot,\,\cdot\, \rangle_x$
represents the canonical pairing between the tangent and cotangent
space.
If $H$ is a Tonelli Hamiltonian, one can easily prove that $L$ is
finite everywhere, of class $C^k$, superlinear and strictly
convex in the fiber in the  $C^2$ sense (\ie $L$ is a 
{\it Tonelli Lagrangian}) and the 
associated Euler-Lagrange flow $\Phi^L_t$ of $L$ is conjugated to the 
Hamiltonian flow $\Phi^H_t$ of $H$ via
the {\it Legendre transform}:
\bea \cL:\; {\rm T}M &\longrightarrow & {\rm T}^*M \nonumber \\
(x,v) &\longmapsto & \left(x,\,\frac{\dpr L }{\dpr v}(x,v)
\right).\label{2.1} \eea
From now on we shall fix $H$ and denote by $L$ its conjugated Lagrangian and,
when referring to an ``invariant'' measure or set, 
we shall understand ``invariant with respect to 
the Hamiltonian flow generated by $H$'' or ``with respect to the 
Euler-Lagrange flow generated by $L$''. 

Given an invariant probability measure $\m$ on ${\rm T}M$ one can 
associate to it an element $\r(\m)$ of the homology group $\rH_1(M;\RRR)$, 
known as {\it rotation vector} or {\it Schwartzman asymptotic cycle},
which generalizes the notion of rotation vector given by 
Poincar\'e and describes how, asymptotically, a 
$\m$-average orbit winds around ${\rm T}M$. 
 In fact, it is easy to show \cite{Mather91} that since $\m$ is invariant by 
the Euler-Lagrangian flow $\Phi^L_t$, if $\eta=df$ is an exact form, then 
$\int{\langle{df},v\rangle} d\m =0$. Therefore, one can define a linear 
functional
\beqano
{\rm H}^1(M;\R) &\longrightarrow& \R \\
c &\longmapsto& \int_{{\rm T}M} \langle{\eta},v\rangle d\m\,,
\eeqano
where $\eta$ is any closed $1$-form on $M$, whose cohomology class is $c$. By 
duality, there exists $\rho (\m) \in {\rm H}_1(M;\R)$ such that
$$
\int_{{\rm T}M} \langle{\eta},v\rangle d\m = \langle c,\rho(\m) \rangle
\qquad \forall\,c\in {\rm H}^1(M;\R).$$ 
$\rho(\m)$ is what we call the {\it rotation vector} of $\m$ and 
it coincides with the Schwartzman asymptotic cycle of $\mu$. 
See \ref{app1} and \cite{Mather91} for more details.

This allows us to define the homology class of certain invariant Lagrangian 
graphs. 

\begin{Def}\rm
A Lagrangian graph $\L$ is called {\it Schwartzman uniquely ergodic} if all 
invariant 
measures supported on $\L$ have the same rotation vector $\r$, which 
will be called {\it homology class of} $\L$. Moreover, if there exists
an invariant measure with full support, $\L$ will be called {\it Schwarztman 
strictly ergodic}.
\end{Def}

We are now ready to state our main result.\\
\\
{\bf{Main Result.}}
{\it For any given $\r\in \rH_1(M;\RRR)$, there exists at most one 
Schwarzman strictly ergodic invariant Lagrangian graph with homology class 
$\r$} [Theorem \ref{uniqschw}, Section \ref{sec3}]\\

For sake of completeness, we also recall the following well-known 
result, which is a corollary of the results in \cite{Mather91} (see Section 
\ref{sec3} below
for a proof).\\
\\
{\bf{Well-known Result.}} {\it For any given $c\in \rH^1(M;\RRR)$, there 
exists at most one invariant Lagrangian
graph $\Lambda$ with cohomology class $c$, carrying an invariant measure 
whose support is the whole of $\Lambda$.}  [Theorem \ref{uniqcom}, 
Section \ref{sec3}]\\

If $M=\TTT^n$, it is natural to ask for the implications of our result 
for KAM theory. In this case, the homology group $\rH_1(\TTT^n;\RRR)$  
is canonically identified with $\R^n$, and  the invariant manifolds of 
interest are the so-called KAM tori, defined as follows.  

\begin{Def}\label{kamtorus}\rm
$\TT\subset\TTT^n\times\RRR^n$ is a (maximal) KAM torus with 
rotation vector $\r$ if:\\
\0i) $\TT\subset\TTT^n\times\RRR^n$ is a Lipschitz  graph over $\TTT^n$;\\
\0ii) $\TT$ is invariant under the Hamiltonian flow $\Phi_t^H$ 
generated by $H$;\\
\0iii) the Hamiltonian flow on $\TT$ is conjugated to a uniform 
rotation on $\TTT^n$; \ie there exists a diffeomorphism 
$\f: \TTT^n\to \TT$
such that $\f^{-1}\circ\Phi_t^H\circ\f= R_\r^t$, $\forall t\in\RRR$,
where $R_\r^t: x \mapsto x+\r t$.
\end{Def}

The celebrated KAM Theorem, whose statement will be recalled in Section 
\ref{sec4}, gives sufficient conditions on $H$ and on the rotation
vector $\r$, allowing one to construct a KAM torus with rotation vector $\r$
and prescribed regularity (depending on the regularity class of $H$). Its 
proof is constructive and the invariant torus one manages to construct is 
locally unique (in a sense that will be clarified in Section \ref{sec4}). 
In spite of the long history and the huge literature dedicated to the KAM 
theorem, the issue of global uniqueness of such tori is still object of some
debate and study, see for instance \cite{BroerTakens}. Our main result 
settles such question in the case of Tonelli Hamiltonians.

\begin{Cor}[Global uniqueness of KAM tori]\label{corollary1}
Every Tonelli Hamiltonian $H$ on ${\rm T}^*\TTT^n$ possesses at most one 
Lagrangian KAM torus for any given rotation vector $\r$. In particular, 
if $H$ and $\r$ satisfy the assumptions of the KAM Theorem, then there exists 
one and only one KAM torus with rotation vector $\r$.
\end{Cor}

The property of being Lagrangian plays a crucial role. As it was observed by 
Herman \cite{Herman}, when $\r$ is
 {\it rationally independent}, \ie 
$\media{\r,\n}\neq 0$, $\forall \n\in\ZZZ^n\setminus\{0\}$, as assumed in the 
KAM theorem, every KAM torus with frequency $\r$ 
is automatically Lagrangian. On the other hand, the existence of
Lagrangian KAM tori with rationally dependent frequency is not typical. 
In some cases, a variant of the classical KAM algorithm allows one 
to construct {\it resonant} invariant tori with a given rational 
rotation vector $\r$, also known as lower dimensional tori 
\cite{GG1,GG2}. However, typically they do not foliate any Lagrangian 
submanifold. Therefore, 
the question of uniqueness of resonant tori is more subtle 
and, to our knowledge, apart from a few partial results \cite{CGGG}, 
it is still open.

In Section \ref{sec4} we shall extend Corollary \ref{corollary1} 
to generic invariant 
tori contained in the $C^0$-closure of the set of KAM tori. 
As remarked by Herman \cite{HermanTori}, 
generically this set is much larger than the
set of KAM tori, and typically the flow on such invariant manifolds is not 
conjugated to a rotation. See Section \ref{sec4} for a more detailed discussion
of these issues. 

%%%%%%%%%%%%%%%%%%%%%%%%%%%%%%%%%%%%%%%%%%%%%%%%%%%%%%%%%%%%%%%%%%%%%%%%%%%%
%%%%%%%%%%%%%%%%%%%%%%%%%%%%%%%%%%%%SECTION 3%%%%%%%%%%%%%%%%%%%%%%%%%%%%%%%
%%%%%%%%%%%%%%%%%%%%%%%%%%%%%%%%%%%%%%%%%%%%%%%%%%%%%%%%%%%%%%%%%%%%%%%%%%%%

\section{Minimizing measures and Lagrangian graphs.}
\label{sec3}

In this section we shall prove the main results announced in Section 
\ref{sec1}. Our proof will be based on Mather's variational approach to the 
study of Lagrangian systems, which is concerned with the study of 
action minimizing
invariant probability measures (also called {\it Mather's measures}) and 
action minimizing orbits of the Euler-Lagrange flow. 
In particular, the keystone of such an approach consists in studying a family 
of modified Tonelli Lagrangians given by 
$L_{\eta}(x,v)=L(x,v)-\langle \eta(x),v\rangle$, where $\eta$ is a closed 
$1$-form on $M$. These Lagrangians, in fact, have all the same Euler-Lagrange 
flow as $L$, but different action minimizing orbits/measures, according to the 
cohomology class of $\eta$. In this way, Mather defined for each cohomology 
class $c \in{\rm H}^1(M;\R)$ two compact invariant subsets of ${\rm T}M$:
\begin{itemize}
\item $\widetilde{\cM}_c$, 
the {\it Mather set of cohomology class $c$}, given 
by the union of the supports of all invariant probability 
measures that minimize the action of $L_{\eta}$ ($c$-{\it action minizimizing 
measures} or {\it Mather's measures of cohomology class} $c$),
where $\eta$ is any closed $1$-form on $M$ with cohomology class $[\eta]=c$; 
\item $\widetilde{\cA}_c$, union of  all {\it regular} 
global minimizers of the action of 
$L_{\eta}$ (or $c$-{\it regular minimizers}); see \cite{Mather93, Fathibook} 
for a precise definition. It is convenient to recall here that, if: 
(i) $\a(c)$ is Mather's $\a$-function, i.e., the maximum over the invariant 
probability measures $\m$ of $-\int L_\eta d\m$, with $[\eta]=c$, and (ii)
$\cS_{\eta}$ is the set of critical subsolutions of $H(x,\eta+p)$, i.e., 
the set of locally Lipschitz functions $u: M \longrightarrow \R$
such that $H_{\eta}(x,d_xu)\leq \a(c)$ for almost every $x\in M$, then 
$\widetilde{\cA}_c={\cal L}^{-1}(\cA^*_c)$, with 
\be {\cA}_c^*=\bigcap_{u\in \cS_{\eta}} \left\{(x,\eta_x + d_xu):\; 
u\;\mbox{is 
differentiable at}\; x\right\} \subset {\rm T}^*M\,,\label{2.7}\ee
see \cite{Fathibook, FathiSiconolfi}.
\end{itemize} 
One can show that $\widetilde{\cM}_c \subseteq \widetilde{\cA}_c$ and, as proved by Carneiro in \cite{Carneiro}, that they are both contained in the energy level
$\tilde{\cal E}_c = \{(x,v)\in {\rm T}M:\, H\circ\cL(x,v) = \a(c)\}$. Moreover, 
one of the most important features of these sets is that they are graphs over 
$M$ (Mather's graph theorem \cite{Mather91}); namely, if $\pi: {\rm T}M 
\rightarrow M$ denotes the canonical projection, then $\pi|{\widetilde{\cA}_c}$
is injective and its inverse $\big(\pi|\widetilde{\cA}_c\big)^{-1}: \cA_c 
\longrightarrow \widetilde{\cA}_c$ is Lipschitz. This is the multidimensional 
analogue of Birkhoff's theorem for twist maps \cite{Birkhoff}.
\\
Analogously, for any rotation vector $h \in{\rm H}_1(M;\R)$, one can define 
another compact invariant subset of ${\rm T}M$:
\begin{itemize}
\item $\widetilde{\cM}^h$, 
the {\it Mather set of homology class $h$}, given by
the union of the supports of all invariant probability measures 
with rotation vector $h$ that minimize the action of $L$ ({\it Mather's 
measures of homology class} $h$).
\end{itemize}
One can show that $\widetilde{\cM}^h$ also enjoys the graph property. 
See Lemma \ref{lemmetto} for the relation between $\widetilde{\cM}^h$ and 
$\widetilde{\cM}_c$.
For more details on Mather's theory for Lagrangian systems, we refer the reader
to \cite{ContrerasIturriaga,Fathibook,Mather91}.

Let us start now  by proving some action minimizing properties of probability 
measures supported on Lagrangian graphs. \\
Given a Lagrangian graph $\Lambda$ with Liouville class $c$, we shall say that
$\Lambda$ is $c$-{\it subcritical}, or simply {\it subcritical}, if 
$\Lambda \subset \{(x,p)
\in {\rm T}^*M:\; H(x,p)\leq \a(c)\}$, where $\a$ is Mather's $\a$-function
(see the lines preceding Eq.(\ref{2.7}) for a definition).
The interest in such graphs comes from the fact that this is the smallest 
energy sub-level of $H$, containing Lagrangian graphs of cohomology $c$ (see 
\cite{Fathibook}). \\ Given a subcritical Lagrangian graph $\Lambda$ with 
Liouville class $c$, 
we shall call $\Lambda_{crit}=\{(x,p)\in \Lambda:\; H(x,p) = \a(c)\}$ its 
{\it critical} part. The key result we need to prove is the following 
characterization of minimizing measures.

\begin{Lem} \label{lemmadellacontesa}
Let $\m$ be an invariant probability measure on ${\rm T}M$ and $\m^*=\cL_*\m$ 
its push-forward to ${\rm T}^*M$, via the Legendre transform $\cL$. Then, 
$\m$ is a Mather's measure if and only if ${\rm supp}\,\m^*$ is contained 
in the critical part of a subcritical Lagrangian graph. In particular, 
any invariant probability measure $\m^*$ on ${\rm T}^*M$, whose support
is contained in an invariant Lagrangian graph with Liouville class $c$,
is the image, via the Legendre transform, of a $c$-action minimizing measure.
\end{Lem}

\begin{Proof}

({\it i})
If $\m$ is a Mather's measure of cohomology class $c$, then the support of 
$\m^*$ is contained in $\cL(\widetilde\cM_c)\subseteq \cL(\widetilde\cA_c)$. 
Since $\cL(\widetilde\cA_c)$ can be obtained as the intersection of all 
$c$-subcritical Lagrangian graphs \cite{bernardc11, Fathibook},
then ${\rm supp}\,\m^*$ is contained in at least one $c$-subcritical Lagrangian
graph $\Lambda$.
In particular ${\rm supp}\,\m^*$ is contained in the critical part of 
$\Lambda$, 
simply because, %from a result by Carneiro \cite{Carneiro}, 
as recalled above,
$\cL(\widetilde\cM_c)$ is contained in the energy level 
${\cal E}^*_c = \{(x,p)\in {\rm T}^*M:\; H(x,p) = \a(c)\}$.

({\it ii}) Let us fix $\eta$ to be a smooth closed $1$-form with $[\eta]=c$,
and let us assume that ${\rm supp}\,\m^*$ 
is contained in the critical part of the $c$-subcritical Lagrangian graph
$\Lambda=\{(x,\eta(x)+du(x))\,,\, x\in M\}$, where $u:M\to\RRR$ is $C^{1,1}$.
In order to show that $\m$ is a $c$-action 
minimizing measure, it is enough to show that 
any orbit $\g$ in ${\rm supp}\,\m$ is a 
$c$-minimizer (i.e., for every finite time-interval $[a,b]$, $\{\g(t)\}_{t\in
[a,b]}$ minimizes the action of $L_\eta$, $[\eta]=c$, 
among the curves with the same end-points $\g(a)$ and $\g(b)$), 
see \cite{Fathibook, ManeI}. 
To this purpose, let us consider
$(x,v)\in{\rm supp}\,\m$ and let $\g(t)\=\p(\Phi_t(x,v))$, where $\Phi_t$
is the Euler-Lagrange flow and $\p$ the canonical projection on $M$. Given any 
interval $[a,b]\subset\RRR$, let us consider the difference 
$u(\g(b))-u(\g(a))$ and rewrite it as:
\be u(\g(b))-u(\g(a))=\int_a^b d_{\g(s)} u(\g(s)) \dot\g(s)\, ds=
\int_a^b \big[L_\eta(\g(s),\dot\g(s))+H_\eta(\g(s),d_{\g(s)} u) \big] ds\,,
\label{4.12}\ee
where the second equality follows from the definition of the Hamiltonian as 
the 
Legendre-Fenchel transform of the Lagrangian and the fact that $\g(s)$ is an 
orbit
of the Euler-Lagrange flow. Note that along the orbit 
$H_\eta(\g(s),d_{\g(s)} u)=\a(c)$, because ${\rm supp}\,\m$ is invariant and
${\rm supp}\,\m^*$ is in the critical part of $\Lambda$. Then
\be \int_a^b L_\eta(\g(s),\dot\g(s)) ds = u(\g(b))-u(\g(a)) -\a(c) (b-a)\;.\ee
On the other hand, any other curve $\g_1:[a,b]\to M$ such that $\g_1(a)=\g(a)$
and $\g_1(b)=\g(b)$ satisfies:
\be u(\g(b))-u(\g(a))=\int_a^b d_{\g_1(s)} u(\g_1(s)) \dot\g(s)\, ds\le
\int_a^b \big[L_\eta(\g_1(s),\dot\g_1(s))+H_\eta(\g_1(s),d_{\g_1(s)} u) \big] 
ds\label{4.14}\ee
where the second inequality follows again by the duality between Hamiltonian 
and Lagrangian. Note that now $H_\eta(\g_1(s),d_{\g_1(s)} u)\le \a(c)$, because
$\Lambda=\{(x,\eta(x)+du(x))\}$ is subcritical. Then
\be \int_a^b L_\eta(\g_1(s),\dot\g_1(s)) ds \ge u(\g(b))-u(\g(a)) -\a(c) (b-a)
\ee
and this proves that $\g$ is a $c$-minimizer. 

Let us finally observe that the Hamilton function on any invariant Lagrangian 
graph $\Lambda=\{(x,\eta+du)\}$ is a constant: $H(x,\eta+du)=k$ 
(the classical proof easily extends to the case of Lipschitz Lagrangian graphs, see for instance \cite{SorrIntegrab}). Then $u$ is a 
classical solution of the Hamilton-Jacobi equation (corresponding to the 
cohomology class $c$). 
As showed in \cite{Fathibook, LPV}, there is only one possible value of $k$ for
which such solutions can exist, namely $k=\a(c)$, and this shows that $\Lambda$
coincides with its critical part. By the result proved in item ({\it ii}), 
if $\m^*$ is supported on $\Lambda$, then $\m$ is a $c$-action minimizing 
measure and this proves the last claim in the statement of the Lemma.
\end{Proof}

We can now prove the well-known uniqueness result  for 
Lagrangian graphs supporting invariant measures of full support, in a fixed 
cohomology class, stated in Section \ref{sec1}.

\begin{Teo}\label{uniqcom}
\it If $\Lambda \subset {\rm T}^*M$  is a Lagrangian graph on which the 
Hamiltonian dynamics admits an 
invariant measure $\m^*$ with full support, then 
$\Lambda = \cL \big(\widetilde{\cM}_c\big)= \cL \big(\widetilde{\cA}_c\big)$, 
where $c$ is the cohomology class of $\Lambda$.\
Therefore, if $\Lambda_1$ and $\Lambda_2$ are two Lagrangian graphs as 
above, with the same cohomology class, then $\Lambda_1=\Lambda_2$.
\end{Teo}

This theorem can be also obtained as a corollary of the results in \cite[Appendix 2]{Mather91}.
Our proof is essentially a rehash of the same ideas, using a different point of view. 
In fact Weierstrass method, used in \cite{Mather91} to show that orbits on a KAM torus are action minimizing, or the use of Hamilton-Jacobi equation are essentially two sides of the same coin.

\begin{Proof} 
By Lemma \ref{lemmadellacontesa},  the measure $\mu=\cL^{-1}_*\m^*$ is 
$c$-action minimizing. This means 
that $\LL^{-1}(\Lambda)={\rm supp}\,\m\subseteq \widetilde\cM_c\subseteq
\widetilde\cA_c$. 
Note however that
$\widetilde\cM_c$ and $\widetilde\cA_c$
are graphs and, since ${\rm supp}\,\m$ is a graph over the 
whole $M$, it follows that
$$ \LL^{-1}(\Lambda)={\rm supp}\,\m = \widetilde\cM_c = \widetilde\cA_c\,.$$
\end{Proof}

One can deduce something more from the above proof. Recall the definition of 
Mather's $\b$-function: $\b(h)$ is the minimum of the average
action $\int L d\m$ over the invariant measures $\m$ with a given
rotation vector $h$ \cite{Mather91}. This function is convex and its {\it 
conjugated} function (given by Fenchel's duality) coincides with the 
$\a$-function: $\a(c)=\max_{h\in {\rm H}_1(M;\R)} 
(\langle c,h\rangle-\b(h))$. 

\begin{Teo}\label{diffbeta}
If $\Lambda$ and $\mu$ are as in Theorem \ref{uniqcom} and $\r$ is 
the rotation vector of 
$\m=\cL^{-1}\m^*$, then $\Lambda=\cL \big(\widetilde{\cM}^\r\big)$. 
Therefore,  if $\Lambda_1$ and $\Lambda_2$ are two Lagrangian graphs 
supporting measures of full support and the same rotation vector $\rho$, 
then $\Lambda_1 = \Lambda_2$. Moreover, Mather's 
$\beta$-function is differentiable at $\r$ with $\partial \beta(\rho) = c$,
where $c$ is the cohomology class of $\L$.\\
\end{Teo}

\begin{Proof}
The first claim follows from the fact that $\widetilde{\cM}^\rho$ is a graph 
over $M$ and that by definition 
$\widetilde{\cM}^\rho \supseteq {\rm supp}\,\m = \LL^{-1}(\Lambda).$
As far as the differentiability of $\beta$ at $\rho$ is concerned,  suppose 
that $c'\in {\rm H}^1(M;\R)$ is a subderivative of $\beta$ at $\rho$.
Then, using Fenchel's duality and the fact that $\a$ and $\b$ are 
conjugated, $\b(\r) = \langle c', \rho \rangle - \a(c')$ and this implies that 
each Mather's measure $\m$ with rotation vector $\r$ is also  $c'$-action 
minimizing; in fact:
$$
\int_{{\rm T}M} (L(x,v)-\langle \eta'(x),v\rangle )\, d\m = \int_{{\rm T}M} 
L(x,v)\,d\m - \int_{{\rm T}M} \langle \eta'(x),v \rangle \, d\m = 
\b(\rho) - \langle c',\rho\rangle = -\a(c')\,, 
$$
where $\eta'$ is a closed $1$-form of cohomology $c'$.
As a result, $\widetilde{\cM}^{\rho} = \cL^{-1}\left(\Lambda\right) \subseteq 
\widetilde{\cM}_{c'}$.
The graph property of $\widetilde{\cM}_{c'}$ and of $\widetilde{\cA}_{c'}$
implies that  $\widetilde{\cA}_{c'} = \widetilde{\cM}_{c'} = \cL^{-1}
\left( \Lambda \right)$ 
and therefore $\cL \big(\widetilde{\cA}_{c'}\big) = \Lambda $. Since the 
cohomology class of $\L$ is $c$, it follows that $c'=c$. 
\end{Proof}

We are now in the position of proving the main uniqueness result in a 
homology class,  stated in Section \ref{sec1}.

\begin{Teo}\label{uniqschw}
Let $\Lambda$ be a Schwartzman strictly ergodic invariant Lagrangian graph with
homology class $\rho$. The following properties are satisfied:
\begin{itemize}
\item[{\rm (i)}] if $\Lambda \cap \cL \big(\widetilde{\cA}_c\big) \neq 
\emptyset$, then $\Lambda=\cL \big(\widetilde{\cA}_c\big)$ and $c=c_{\Lambda}$,
where $c_{\Lambda}$ is the cohomology class of $\Lambda$.
\item[{\rm (ii)}] the Mather function $\alpha $ is differentiable at 
$c_{\Lambda}$ 
and $\partial \alpha (c_{\Lambda}) = \rho$.
\end{itemize}
Therefore,
\begin{itemize}
\item[{\rm (iii)}] any invariant Lagrangian graph that carries a measure with 
rotation vector $\rho$ is equal to the graph $\Lambda$;
\item[{\rm (iv)}] any invariant Lagrangian graph is either disjoint from 
$\Lambda$ or equal to $\Lambda$.
\end{itemize}
\end{Teo}

We shall need the following Lemma.

\begin{Lem}\label{lemmetto} Let $\rho,c$ be respectively an arbitrary homology 
class in $\rH_1(M;\R)$ and an arbitrary cohomology class in $ \rH^1(M;\R)$. We 
have
$$
{\rm (1)}\; \widetilde{\cM}^{\rho} \cap  \widetilde{\cM_c} \neq \emptyset 
\quad \Longleftrightarrow \quad
{\rm (2)}\; \widetilde{\cM}^{\rho} \subseteq \widetilde{\cM_c} \quad 
\Longleftrightarrow \quad {\rm (3)}\;\rho \in \partial \a(c)\,.
$$
\end{Lem}

\noindent{\bf Proof of Lemma \ref{lemmetto}.} The implication
${\rm (2) }\Longrightarrow {\rm (1)}$ is trivial. Let us prove that
 ${\rm  (1)} \Longrightarrow {\rm (3)}$.
If $\widetilde{\cM}^{\rho} \cap  \widetilde{\cM_c}\neq 0$, then there exists a 
$c$-action minimizing invariant measure $\m$ with rotation vector $\rho$. 
Let $\eta$ be a closed $1$-form with $[\eta]=c$; 
from the definition of $\a$ and $\b$:
\beqano
-\a(c) = \int_{{\rm T}M} (L(x,v)-\langle{\eta(x),v}\rangle)\,d\m = 
\int_{{\rm T}M} L(x,v)\,d\m - 
\langle c,\rho\rangle = \b(\rho) - \langle c,\rho\rangle\,;
\eeqano
since $\b$ and $\a$ are convex conjugated,  then $\rho$ is a subderivative of 
$\a$ at $c$.\\
Finally, in order to show $(3) \Longrightarrow (2)$, let us prove that any 
Mather's measure with rotation vector $\rho$ is $c$-action minimizing. 
In fact, if $\rho \in \partial \a(c)$ then $\a(c)= \langle c,\rho\rangle - 
\b(\rho)$; therefore for any $\m$ is a Mather's measure with rotation vector 
${\rho}$ and $\eta$ as above:
$$
-\a(c) = \b(\rho) - \langle c,\rho\rangle =  \int_{{\rm T}M}(L(x,v)-
\langle{\eta(x),v}\rangle)\, d\m.  
$$
This proves that $\m$ is $c$-action minimizing and concludes the proof.
\qed\\

\noindent{\bf Proof of Theorem \ref{uniqschw}} 
({\it i}) From Theorem \ref{uniqcom}, it follows that $\cL^{-1} \big(\Lambda 
\big)=\widetilde{\cA}_{c_{\Lambda}}$.
Let us show that it does not intersect any other Aubry set.
Suppose by contradiction that $\cL^{-1} \big(\Lambda \big)$ intersects another 
Aubry set  
$\widetilde{\cA}_{c}$.
By Theorem \ref{diffbeta}, $\cL^{-1} \big(\Lambda \big) = 
\widetilde{\cM}^{\rho}$, then 
$\widetilde{\cM}^{\rho} \cap \widetilde{\cA}_c \neq \emptyset$ and, because 
of lemma \ref{lemmadellacontesa}, lemma \ref{lemmetto} and the graph property 
of $\widetilde{\cA}_c$, we can 
conclude that $\widetilde{\cA}_c = \cL^{-1} \big(\Lambda \big)$. The same 
argument used in the proof of Theorem \ref{diffbeta} allows us to conclude that
$c=c_{\Lambda}$.\\
({\it ii}) Suppose that $ h \in \partial \a(c_{\Lambda}) $. The previous lemma 
implies that $\widetilde{\cM}^{h} \subseteq \Lambda $; the Schwartzman unique 
ergodicity property of $\Lambda$ implies $h=\rho$. Therefore
$\alpha $ is differentiable at $c_{\Lambda}$ and $\partial\alpha 
(c_{\Lambda})=\rho$.

To prove ({\it iv}), let $\Lambda_1$ be an invariant Lagrangian graph, and call
$c_1$ its cohomology class.
If the compact invariant set $\Lambda\cap\Lambda_1$ is not empty, 
then we can find a probability measure $\mu^*$ invariant under the flow
and whose support is contained in this intersection. Since $\mu^*$ is contained
in the Lagrangian graph $\Lambda_1$, by Lemma \ref{lemmadellacontesa}, 
it is $c_1$-action minimizing. Hence, the support of $\mu^*$ is contained in 
$\cL \big(\widetilde{\cA}_{c_1}\big)$.
This shows that the intersection $\Lambda \cap \cL \big(\widetilde{\cA}_{c_1}
\big)$ contains the support of $\mu^*$ and is therefore not empty. 
By ({\it i}), $\Lambda=\cL \big(\widetilde{\cA}_{c_1}\big)$. Moreover, note
that $\cL \big(\widetilde{\cA}_{c_1}\big)\subseteq
\L_1$, because $\L_1={\rm Graph}(\eta_1+du_1)$, with $[\eta_1]=c_1$ and $u_1$ a
classical 
solution to the Hamilton-Jacobi equation (see \cite{Fathibook}).
Therefore, $\Lambda=\Lambda_1$, since they are both graphs over $M$.

To prove ({\it iii}), consider an invariant Lagrangian graph $\Lambda_1$, 
with cohomology class $c_1$, which carries an invariant measure 
$\mu^*$ whose rotation vector is $\rho$. By Lemma \ref{lemmadellacontesa},
the measure $\mu^*$  is $c_1$-minimizing. Therefore, we have 
$\widetilde{\cM}^\rho\cap \widetilde\cM_{c_1}\neq \emptyset$. 
By Lemma \ref{lemmetto}, it 
follows that  ${\cal L}^{-1}(\Lambda)=\widetilde{\cM}^\rho\subseteq 
\widetilde\cM_{c_1}\subseteq 
\widetilde\cA_{c_1}\subseteq {\cal L}^{-1}(\Lambda_1)$. Again, this forces the 
equality $\Lambda=\Lambda_1$ by the graph property.\qed
\\

Finally, observe that using Lemma \ref{lemmetto}, one can also deduce the 
following property.

\begin{Cor} Mather's $\alpha$-function is differentiable at $c$ if and only if 
the restriction of the Euler-Lagrange flow to $\widetilde{\cA}_c$ 
is Schwartzman uniquely ergodic, \ie if and only if all invariant measures 
supported on $\widetilde{\cA}_c$ have the same rotation vector.
\end{Cor}

%%%%%%%%%%%%%%%%%%%%%%%%%%%%%%%%%%%%%%%%%%%%%%%%%%%%%%%%%%%%%%%%%%%%%%%%%%%
%%%%%%%%%%%%%%%%%%%%%%%% SECTION 4%%%%%%%%%%%%%%%%%%%%%%%%%%%%%%%%%%%%%%%%%
%%%%%%%%%%%%%%%%%%%%%%%%%%%%%%%%%%%%%%%%%%%%%%%%%%%%%%%%%%%%%%%%%%%%%%%%%%%

\section{Global uniqueness of KAM tori}\label{sec4}

%%%%%%%%%%%%%%%%%%%%%%%%%%%%%%%%%%%%%%%%%%%%%%%%%%%%%%%%%%%%%%%%%%%%%%%%%%%

In this section we motivate more precisely the problem of uniqueness of KAM 
tori and prove Corollary \ref{corollary1}.
We also show how to generalize Corollary \ref{corollary1} to cover the case of
invariant tori belonging to the closure of the set of KAM tori. 

KAM theory concerns the study of existence of KAM tori (see Definition 
\ref{kamtorus}) in quasi-integrable Hamiltonian systems of the form
$H(x,p)=H_0(p)+\e f(x,p)$, where: $(x,p)$ are local coordinates on $\TTT^n
\times\RRR^n$, $\e$ is a ``small''  parameter and $f(x,p)$ a smooth function.
If $\e=0$ the system is integrable, in the sense that the dynamics
can be explictly solved: in particular each torus $\TTT^n\times\{p_0\}$
is invariant and the motion on it corresponds to a rotation 
with frequency $\r(p_0)=\frac{\partial H_0}{\partial p}(p_0)$. The question 
addressed by KAM theory is whether this foliation of phase space into
invariant tori, on which the motion is (quasi-)periodic, persists even 
if $\e\neq 0$. 
In 1954 Kolmogorov stated (and Arnol'd \cite{Arnold} and Moser \cite{Moser} 
proved it later in different 
contexts) that, in spite of the generic disappearence of the invariant 
submanifolds filled by periodic orbits, pointed out by Poincar\'e, 
for small $\e$ it is always 
possible to find KAM tori corresponding to ``{\it strongly non-resonant}'',
\ie {\it Diophantine}, rotation vectors. 
The celebrated KAM Theorem (in one of its several versions) not only shows the 
existence of such tori, but also provides an explicit method to construct them.
\\

\noindent{\bf Theorem} (Kolmogorov--Arnol'd--Moser) \cite{Salamon}. 
{\it Let $n\ge 2$, $\t>n-1$, $C>0$, $\ell>2\t+2$, $M>0$
and $r>0$ be given. Let $B_r\in\RRR^n$ be the open ball of radius $r$ 
centered at the origin. Let $H\in C^\ell(\TTT^n\times B_r)$ be of the form 
\be H(x,p)=H_0(p)+\e f(x,p)\label{1.2}\ee
with $|H_0|_{C^\ell}\le M$, $|f|_{C^\ell}\le M$, $\left|\frac{\partial^2 H_0}
{\partial p^2}\right|\ge M^{-1}$ and
$\r= \frac{\partial H_0}{\partial p}(0)\in \DD(C,\t)$. Then, for any $s<\ell-
2\t-1$, 
there exists $\e_0>0$ such that for
any $\e\le\e_0$ the Hamiltonian {\rm(\ref{1.2})} admits a $C^{s,s+\t}$ 
KAM torus with rotation vector $\r$, \ie a $C^{s+\t}$ invariant torus such that
the Hamiltonian flow on it is $C^s$-conjugated with a rotation with 
frequency $\r$.} \\

The invariant torus constructed in the proof of the KAM Theorem is locally
unique, in the sense that for any prescribed (and admissible) $s$ 
there is at most one $C^{s,s+\t}$ KAM torus with rotation 
vector $\r$ within a $C^s$-distance $\d(n,s,C,\t)$ to the one constructed 
in the proof of the KAM Theorem, see \cite{BroerTakens, 
salamon-zehnder, Salamon}. Note that the $C^s$-distance 
$\d$ within which one can prove uniqueness of the KAM torus in a prescribed
regularity class depends both on the irrationality properties of $\r$ and
on the regularity class $s$ itself. It is then {\it a priori} possible that 
even for small $\e$ there exist different KAM tori {\it within a prescribed 
$C^1$-distance} from the one constructed in the proof of the Theorem, possibly 
less regular than that torus. Quite surprisingly, even in the analytic case,
we are not aware of any proof of ``global'' uniqueness of the 
invariant analytic KAM torus with rotation vector $\r$ (of course in the 
analytic case the analytic torus one manages to construct is unique within
the class of analytic tori -- however nothing {\it a priori} guarantees that 
less regular invariant tori with the same rotation vector exist).
\\

Our result, in the form stated 
in Corollary \ref{corollary1}, settles the question and shows that,
at least in the case of Tonelli Hamiltonians, it is not possible
to have two different KAM tori with the same rotation vector. Note that
the assumption of strict convexity of the Hamiltonian is necessary
to exclude trivial sources of non-uniqueness: for instance, in 
the context of quasi-integrable Hamiltonians, 
global uniqueness could be lost simply because the unperturbed Hamiltonian 
induces a  map $p\mapsto \partial_p H_0(p)$ from actions to frequencies that is
not one to one. Let us also remark that, 
apparently, the Hamiltonian considered in KAM Theorem is not a Tonelli
Hamiltonian, since the latter, by definition, is defined globally on the whole 
$\TTT^n\times\RRR^n$. However any $C^\ell$ strictly convex Hamiltonian defined
on $\TTT^n\times B_r$ for some $r>0$ can be extended to a global
$C^\ell$ Tonelli Hamiltonian. Then in the statement of the KAM Theorem above 
it is actually enough to assume $H$ to be a $C^\ell$ Tonelli Hamiltonian, 
locally satisfying the (in)equalities listed after (\ref{1.2}).\\

Given the proof of our main results in Section \ref{sec3}, 
the proof of Corollary \ref{corollary1} becomes very simple.\\

\noindent{\bf Proof of Corollary \ref{corollary1}.}
Since the Lagrangian KAM torus $\TT$ admits an invariant measure $\m^*$ of full
support, which is the image via the conjugation $\f$ of the uniform 
measure on $\TTT^n$, then the claims follow from Theorem \ref{diffbeta}. 
Note that for rationally independent rotation vectors,
a classical remark by Herman \cite{Herman} 
implies that $\TT$ is automatically Lagrangian (it is sufficient that the flow 
on it is topologically conjugated to a transitive flow on $\T^n$).\qed\\

An interesting generalization of the result of Corollary \ref{corollary1}
concerns the invariant tori belonging to the $C^0$-closure $\overline\Upsilon$ 
of the set $\Upsilon$ of all Lagrangian KAM tori. Note that, for
quasi-integrable systems, $\Upsilon$ is not empty.
The set $\Upsilon$ can be seen as a subset of ${\rm Lip}(\TTT^n,\RRR^n)$.
This follows from Theorem \ref{diffbeta},
and from Mather's graph theorem and the other results in \cite{Mather91}. 
Moreover, any family of invariant Lagrangians 
graphs on which the function $\alpha$ (or $H$) is bounded gives rise to a 
family of functions in ${\rm Lip}(\TTT^n,\RRR^n)$ with uniformly bounded 
Lipschitz constant. This is because, given $\L$ in such a family and 
denoting by $(\eta+du)$ its graph, for any pair 
of points $x,y\in\TTT^n$ and any smooth curve $\g(t)$ on $\L$ 
connecting $x$ to $y$ with unit speed, we have that $u(x)-u(y)\le
\int_0^{|x-y|}L_\eta(\g(t),\dot\g(t))+\a(c)|x-y|$, where $c=[\eta]$, 
see (\ref{4.14}).
By Ascoli-Arzel\`a theorem, it follows  that $\overline\Upsilon$
is also a subset of ${\rm Lip}(\TTT^n,\RRR^n)$, consisting of  functions whose 
graphs are invariant  (Lipschitz) Lagrangian tori.
Herman \cite{HermanTori}
showed that, for a generic Hamiltonian $H$ close enough to an 
integrable Hamiltonian $H_0$, the dynamics on the generic  tori in  $\overline
\Upsilon$ is not conjugated to a 
rotation. These  ``new'' tori therefore represent the majority, in the sense 
of topology, and hence most invariant tori cannot be obtained by the KAM 
algorithm. 
More precisely, Herman showed that in $\overline\Upsilon$ there exists a dense 
$G_\d$-set (\ie
a dense countable intersection of open sets) of invariant Lagrangian graphs  
on which the dynamics is strictly ergodic and weakly mixing, and for which the 
rotation vector, in the sense of Section \ref{sec1}, is not Diophantine. These 
invariant graphs are therefore not obtained by the KAM theorem, however our 
uniqueness result do still apply to these graphs since strict ergodicity 
implies Schwartzman strict ergodicity.

More generally, given any Tonelli Lagrangian on $\T^n$, we consider the set 
$\tilde \Upsilon$ of invariant Lagrangian graphs on which the dynamics 
of the flow is topologically conjugated to an {\it ergodic} linear flow on 
$\T^n$ 
(of course, far from the canonical integrable Lagrangian the set $\tilde 
\Upsilon$ may 
be empty). The dynamics on anyone of the invariant graphs in $\tilde \Upsilon$ 
is strictly ergodic. Since the set of strictly ergodic flows on a compact set 
is a 
$G_\delta$-set in the $C^0$ topology, see for example \cite[Corollaire 4.5]
{FathiHerman}, 
it follows that there exists a dense $G_\delta$-subset ${\cal G}$ of the $C^0$ 
closure   
of $\tilde \Upsilon$ in  ${\rm Lip}(\TTT^n,\RRR^n)$, such that the dynamics on 
any $\Lambda\in {\cal G}$ is strictly ergodic. 
Therefore we get the following proposition.

\begin{Prop}
There exists a dense $G_\d$-set $\cal G$ in the $C^0$ closure of $\tilde
\Upsilon$ consisting of 
strictly ergodic invariant Lagrangian graphs.  Any $\Lambda\in {\cal G}$ 
satisfies the following properties:
\begin{itemize}
\item[{\rm (i)}] the invariant graph $\Lambda$ has a well-defined rotation 
vector $\rho(\Lambda)$.
\item[{\rm (ii)}]  Any invariant Lagrangian graph that intersects $\Lambda$ 
coincides with $\Lambda$.
\item[{\rm (iii)}] Any Lagrangian invariant graph that carries an invariant 
measure whose rotation is $\rho(\Lambda)$ coincides with $\Lambda$. 
\end{itemize}
\end{Prop}

%%%%%%%%%%%%%%%%%%%%%%%%%%%%%%%%%%%%%%%%%%%%%%%%%%%%%%%%%%%%%%%%%%%%%%%%%%%%

\appendix
\newcommand{\appsection}[1]{\let\oldthesection\thesection
  \renewcommand{\thesection}{Appendix \oldthesection}
  \section{#1}\let\thesection\oldthesection}

\appsection{Schwartzman unique and strict ergodicity}\label{app1}

In this section we prove some results on the structure of the set of 
Schwartzman uniquely/strictly ergodic graphs that we have
introduced in Section \ref{sec1}, and provide some examples.

Let us start from the notion of Schwartzman asymptotic cycle of a flow, 
introduced by Sol Schwartzman in \cite{Schwartzman}, as a first attempt to 
develop an algebraic topological approach to the study of dynamics. This is 
closely related to the concept of rotation vector of a measure, that we have 
introduced in Section \ref{sec1}.
We shall give a different description of the Schwartzman asymptotic cycle of a 
flow using the flux homomorphism in volume preserving and symplectic
geometry (see \cite{Banyaga}[Chapter 3]), from the same perspective as  
\cite{Fathi80}. The definition used below has the technical advantage of 
not relying on the Krylov-Bogolioubov theory of generic orbits in a dynamical 
system,  although a more geometrical definition showing that ``averaged'' 
pieces of long orbits converge almost everywhere in the first homology group 
for any invariant measure is certainly more heuristic and intuitive.\\

Let $X$ be a topological space and $(\phi_t)_{t\in\R}$ a 
continuous flow on $X$. We shall define $\Phi:X\times [0,1]\to X$ by $\Phi(x,t)
=\phi_t(x)$.
Consider a continuous function $f:X\to \T$ and let
$F(f,\Phi):X\times[0,1]\to\T$ be 
$$F(f,\Phi)(x,t)= f(\phi_t(x))-f(x).$$
$F(f,\Phi)$ is continuous and identically $0$ on $X\times\{0\}$, it is 
therefore homotopic to a constant and can be lifted to a continuous map 
$\bar{F}(f,\Phi):X\times[0,1]\to\R$, with $\bar F(f,\Phi)|X\times\{0\}$ 
identically $0$. We define 
$${\cal V}(f,\phi_t)(x):=\bar F(f,\Phi)(x,1).$$
Note that if $f$ is homotopic to $0$ then it can be lifted continuously to 
$\bar f:X\to \R$.
In that case $\bar F(f,\Phi)=\bar f\Phi-\bar f$, and
$${\cal V}(f,\phi_t)(x)=\bar f(\phi_1(x))-\bar f.$$
If $\mu$ is a measure with compact support invariant under the flow $\phi_t$, 
for a continuous 
$f:X\to\T$, we define
${\cal S}(\mu,\phi_t)(f)$, or simply ${\cal S}(\mu)(f)$ when $\phi_t$ is fixed,
by $${\cal S}(\mu)(f)=\int_X{\cal V}(f,\phi_t)(x)\,d\mu(x).$$
If we denote by $[X,\T]$ the set of homotopy classes of continuous 
maps from $X$ to $\T$, it is not difficult to verify that 
${\cal S}(\mu)$ is a  well-defined additive homomorphism from the additive 
group $[X,\T]$ to $\R$.

When $X$ is a good space (like a manifold or a locally finite polyhedron), it 
is well-known that $[X,\T]$ is canonically identified with the first cohomology
group $\rH^1(X;\Z)$.
In that case ${\cal S}(\mu)$ is in $\operatorname{Hom}(\rH^1(X;\Z),\R)$. Since 
the first cohomology group with real coefficients $\rH^1(X;\R)$ is $\rH^1(X;\Z)
\otimes \R$,
we can view ${\cal S}(\mu)$ as an element of the dual $\rH^1(X;\R)^*$ of the 
$\R$-vector space $\rH^1(X;\R)$. When $\rH^1(X;\R)$ is finite-dimensional (for 
instance, when $X$ is a finite polyhedron or a compact manifold)
then $\rH^1(X;\R)^*$ is in fact equal to the first homology group 
$\rH_1(X;\R)$, and therefore ${\cal S}(\mu)$ defines an element of 
$\rH_1(X;\R)$, \ie a 1-cycle. 
This 1-cycle ${\cal S}(\mu)$ is called the {\it Schwartzman asymptotic cycle} 
of $\mu$. \\

Let us now consider the case of a 
$C^1$ flow $\phi_t$ on a manifold $N$. We call $X$ the continuous vector 
field on $N$ 
generating $\phi_t$, \ie
$$\forall x\in N, \quad X(x)={\frac {d\phi_t(x)}{dt}}\Big|_{t=0}.$$
By the flow property $\phi_{t+t'}=\phi_t\circ\phi_{t'}$, this implies 
$$
\forall x\in N,\ \forall t\in \R, \quad \frac {d\phi_t(x)}{dt}=X(\phi_t(x)).
$$

In the case of a manifold $N$, the identification of $[N,\T]$ with 
$\rH^1(N;\Z)$ is best described with the de Rham cohomology. We consider the 
natural map 
$I_N:[N,\T]\to \rH^1(N;\R)$ defined by 
$$I_N([f])=[f^*\theta],$$
where $[f]$ on the left hand side denotes the homotopy class of the 
$C^\infty$ map $f:N\to\T$, and $[f^*\theta]$ on the right hand side  
is the cohomology class of the pullback by $f$ of the closed 1-form on $\T$ 
whose lift to $\R$ is $dt$. Note that any homotopy class in $[N,\T]$ contains 
smooth maps because $C^\infty$ maps are dense in $C^0$ maps (for the Whitney 
topology). Therefore the map $I_N$ is indeed defined on 
the whole of $[N,\T]$.
As it is well-known, this map $I_N$ induces an isomorphism of $[N,\T]$ on 
$\rH^1(N;\Z)\subset \rH^1(N;\R)=\rH^1(N;\Z)\otimes \R$.

Given a $C^\infty$ map $f:N\to \T$, the $C^1$ flow $\phi_t$ on $N$, and 
$x\in N$, we compute ${\cal V}(f,\phi_t)(x)$. If $\gamma_x:[0,1]\to N$ is the 
path
$t\mapsto\phi_t(x)$, 
since $\gamma_x$ is $C^1$, we get
$${\cal V}(f,\phi_t)(x)=\int_{f\circ \gamma_x}\theta=\int_{\gamma_x}f^*\theta.
$$
Moreover, since $\gamma_x(t)=\phi_t(x)$, we have $\dot\gamma_x(t)=
X(\phi_t(x))$. It follows that 
$${\cal V}(f,\phi_t)(x)=\int_0^1(f^*\theta)_{\phi_t(x)}(X[\phi_t(x)])\,dt = 
\int_0^1(i_Xf^*\theta)(\phi_t(x))\,dt,$$
where $i_X$ denotes the interior product of a differential form with $X$.
Therefore if $\mu$ is an invariant measure for $\phi_t$, which we shall assume 
to have a compact support, we obtain
\beqano
{\cal S}(\mu)&=&\int_N\int_0^1(i_Xf^*\theta)(\phi_t(x))\,dtd\mu(x)
= \int_0^1\int_N(i_Xf^*\theta)(\phi_t(x))\,d\mu(x)dt = \\
&=& \int_0^1\int_N(i_Xf^*\theta)(x)\,d\mu(x)dt=\int_N(i_Xf^*\theta)
\,d\mu.
\eeqano
This shows that as an element of $\rH^1(M;\R)^*$, the Schwarztman asymptotic 
cycle 
${\cal S}(\mu)$ is given by
$${\cal S}(\mu)([\eta])=\int_Ni_X\eta\,d\mu.$$\\

We would like now to relate the Schwartzman asymptotic cycles to the rotation 
vectors $\rho(\mu)$ defined in section \ref{sec1} for Euler-Lagrange flows. 
In this case $N=TM$ and $\phi_t$ is an Euler-Lagrange flow
$\phi_t^L$ of some Lagrangian $L$. If we call $X_L$ the vector field generating
$\phi_t^L$, since this flow is obtained from a second order ODE on $M$, we get
$$\forall x\in M, \forall v\in T_xM, \quad T\pi(X_L(x,v))=v,$$
where $T\pi:T(TM)\to TM$ denotes the canonical projection.
Since this projection $\pi$ is a homotopy equivalence, to compute 
${\cal S}(\mu)$ we only need to consider forms of the type $\pi^*\eta$ where 
$\eta$ is a closed 1-form on the base $M$. In this case $(i_{X_L}\pi^*\eta)
(x,v)=\eta_x(T\pi(X_L(x,v))=\eta_x(v)$. Therefore, for any  
probability measure $\mu$ on $TM$ with compact support and invariant under 
$\phi_t^L$,  we obtain
$${\cal S }(\mu)[\pi^*\eta]=\int_{TM}\eta_x(v)\,d\mu(x,v)
=\int_{TM}\langle\eta, v \rangle\,d\mu.$$
This is precisely $\rho(\mu)$ as it was defined above in section \ref{sec1}. 
Note that the only property we have used is the fact that $\phi_t$ is the flow 
of a second 
order ODE on the base $M$.\\

\noindent{\bf Examples}.
\begin{itemize}
\item[-] We can now easily compute Schwartzman asymptotic cycles for linear 
flows on 
$\T^n$. Such a flow is determined by a constant vector field $\alpha\in \R^n$ 
on $\T^n$ (here we use the canonical trivialization of the tangent bundle of 
$\T^n$), the associated flow $R^\alpha_t:
\T^n\to\T^n$ is defined by $R^\alpha_t(x)=x+[t\alpha]$, where $[t\alpha]$ is 
the class in $\T^n=\R^n/\Z^n$ of the vector $t\alpha\in\R^n$. If $\omega$ is a 
1-form with constant coefficients, \ie $\omega=\sum_{i=1}^na_idx_i$, with 
$a_i\in\R$, the interior product
$i_\alpha\omega$ is the constant function $\sum_{i=1}^n\alpha_ia_i$. Therefore,
it follows that ${\cal S}(\mu)=\alpha\in \R^n\equiv \rH_1(\T^n;\R)$. 
\item[-] Suppose that $x$ is a periodic point of $\phi_t$ or period $T>0$. One 
can define an invariant  probability measure $\mu_{x,t_0}$ for $\phi_t$ by
$$\int_Xg(x)\,d\mu_{x,t_0}=\frac1{t_0}\int_0^{t_0}g((\phi_t(x))\,dt,$$
where $g:X\to \R$ is a measurable function. We let the reader verify that
${\cal S}(\mu_{x,t_0})$ is equal in $\rH_1(X;\R)$ to the homology class 
$[\gamma_{x,t_0}]/t_0$, where $\gamma_{x,t_0}$ is the loop $t\mapsto \phi_t(x),
t\in[0,t_0]$.
\item[-] When $x$ is a fixed point of $\phi_t$, then the Dirac mass 
$\delta_x$ at $x$ is invariant under $\phi_t$, and in that case ${\cal S}
(\delta_x)=0$.\\
\end{itemize}
 
Let us now study the behavior of Schwartzman asymptotic cycles under 
semi-conjugacy. 
\begin{Prop}\label{SchSemiConj}
Suppose $\phi^i_t:X_i\to X_i, i=1,2$ are two continuous flows. Suppose also 
that $\psi:X_1\to X_2$ is a continuous semi-conjugation between the flows, \ie
$\psi\circ \phi^1_t=\phi^2_t\circ \psi$, for every $t\in\R$. Given a 
probability measure
$\mu$ with compact support on $X_1$ invariant under $\phi^1_t$, then, 
for every continuous map $f:X_2\to\T$, we have
$${\cal S}(\psi_*\mu,\phi^2_t)([f])={\cal S}(\mu,\phi^1_t)([f\circ \psi]),$$
where  $\psi_*\mu$ is the image of $\mu$ under $\psi$. In particular, if we 
are in the situation where 
$\operatorname{Hom}([X_i,\T])$ $\equiv \rH_1(X_i;\R),
i=1,2$, we obtain
$${\cal S}(\psi_*\mu,\phi^2_t)=H_1(\psi)({\cal S}(\mu,\phi^1_t)).$$
\end{Prop}
\begin{Proof} Notice that $f\psi\phi^1_t(x)-f\psi(x)=f\phi^2_t(\psi(x))-
f(\psi(x))$. Therefore by uniqueness of liftings ${\cal V}(f\psi,\phi^1_t)(x)
={\cal V}(f,\phi^2_t)(\psi(x))$.
An integration with respect to $\mu$ finishes the proof.
\end{Proof}

To simplify things, in the remainder of this appendix, we shall assume that 
$X$ is a compact space, for which we have $[X,\T]=\rH^1(X;\Z)$, and 
$\rH_1(X;\Z)$ is finitely generated. In that case,
the dual space $\rH^1(X;\R)^*$ is $\rH_1(X;\R)$, and for every flow $\phi_t$ 
on $X$ and every probability measure $\mu$ on $X$ invariant under $\phi_t$, 
the Schwartzman asymptotic cycle is an element of the finite dimensional-vector
space $\rH_1(X;\R)$.

\begin{Def}\rm For a flow $\phi_t$ on $X$, we denote by ${\cal S}(\phi_t)$ 
the set of all 
Scwhartzman asymptotic cycles ${\cal S}(\mu)$, where $\mu$ is an arbitrary 
probability measure on $X$ invariant under $\phi_t$.
\end{Def}

Since $X$ is compact, note that for the weak topology the set $\calM(X)$ of 
probability Borel measures on $X$ is compact and convex. It is even metrizable,
since we are assuming $X$ metrizable.
Furthermore the subset ${\calM}(X,\phi_t)\subseteq {\calM}(X)$ of probability 
measures invariant under $\phi_t$ is, as it is well-known, compact convex and 
non empty. Therefore
 ${\cal S}(\phi_t)$ is a compact convex non-empty subset of $\rH_1(X;\R)$.
 For the case of a linear flow $R^\alpha$ on $\T^n$, we have shown above that
${\cal S}(R^\alpha_t)=\{\alpha\}\subset \R^n\equiv \rH_1(\T^n;\R)$.
 
The following corollary is an easy consequence of Proposition 
\ref{SchSemiConj}.
\begin{Cor}\label{SchCon} For $i=1,2$, suppose that $\phi^i_t$ is a continous 
flow on the compact space $X_i$, which satisfies $\operatorname{Hom}([X_i,\T],
\R)\equiv \rH_1(X_i;\R)$.
If $\psi:X_1\to X_2$ is a topological conjugacy between $\ph^1_t$ and 
$\phi^2_t$
(\ie the map $\psi$ is a homeomorphism that satisfies $\psi\phi^1_t=\phi^2_t 
\psi$, for all $t\in \R$), then we have
$${\cal S}(\phi^2_t)=H_1(\psi)[{\cal S}(\phi^1_t)].$$
\end{Cor}  

We denote by ${\mathfrak F}(X)$ the set of continuous flows on $X$. We can 
embed ${\mathfrak F}(X)$ in $C^0(X\times [0,1],X)$ by the map $\phi_t\mapsto 
F^{\phi_t}\in 
C^0(X\times [0,1],X)$, where
$$F^{\phi_t}(x,t)=\phi_t(x).$$
The topology on $C^0(X\times [0,1],X)$ is the compact open (or uniform) 
topology, and
we endow ${\mathfrak F}(X)$ with the topology inherited from the embedding 
given above.

\begin{Lem}The map $\phi_t\mapsto {\cal S}(\phi_t)$ is upper semi-continuous 
on ${\mathfrak F}(X)$. This means that for each open
subset $U\subseteq \rH_1(X;\R)$, the set 
$\{\phi_t\in {\mathfrak F}(X)\mid {\cal S}(\phi_t)\subset U\}$ is open in 
${\mathfrak F}(X)$.
\end{Lem}
\begin{Proof} Since the topology on $C^0(X\times [0,1],X)$ is metrizable, 
if this were not true we could find an open set $U\subset \rH_1(X;\R)$ and a 
sequence  $\phi^n_t$ of continuous flows on $X$ converging uniformly to a flow
$\phi_t$, with ${\cal S}(\phi_t)\subset U$, and 
$ {\cal S}(\phi^n_t)$ is not contained in $ U$. This means that for each $n$ we
can find a probability measure $\mu_n$ on $X$ invariant under $\phi^n_t$ and 
such that
its Schwartzman asymptotic cycle ${\cal S}(\mu_n,\phi^n_t)$ for $\phi_n^t$ is 
not in the open set $U$. Since ${\calM}(X)$ is compact for the weak topology, 
extracting a subsequence if necessary, we can assume that $\mu_n\to\mu$. It is 
not difficult to show that
$\mu$ is invariant under the flow $\phi_t$. We now show that ${\cal S}(\mu_n,
\phi^n_t)\to {\cal S}(\mu,\phi_t)$. This will yield a contradiction and finish 
the proof because ${\cal S}(\mu_n,\phi^n_t)$ is in the closed set $\rH_1(X;\R)
\setminus U$, for every $n$, and ${\cal S}(\mu,\phi_t)\in U$.

To show that the linear maps ${\cal S}(\mu_n,\phi^n_t)\in \rH_1(X;\R)=
\rH^1(X;\R)^*$ converge to the linear map ${\cal S}(\mu,\phi_t)$, it suffices 
to show that 
${\cal S}(\mu_n,\phi^n_t)([f])\to {\cal S}(\mu,\phi_t)([f])$, 
for every $[f]\in[X,\T]=\rH^1(X;\Z)\subset \rH^1(X;\R)=\rH^1(X;\Z)\otimes\R$.
Fix now a continuous map $f:X\to \T$. Denote by $F_n, F:X\times [0,1]\to \T$ 
the maps
defined by 
$$F_n(x,t)=f(\phi^n_t(x))-f(x)\text { and}\; F(x,t)=f(\phi_t(x))-f(x).$$
By the uniform continuity of $f$ on the compact metric space $X$,
the sequence $F_n$ converges uniformly to $F$. Since 
$F_n|X\times \{0\}\equiv 0$, if we call $\tilde F_n:X\times [0,1]\to \R$ 
the lift of $F_n$ such that $\tilde F_n|X\times \{0\}\equiv 0$, then the 
sequence 
$\tilde F_n$ also converges uniformly to $\tilde F$, that is the lift of $F$ 
such that  $\tilde F |X\times \{0\}\equiv 0$. Since the $\mu_n$ are probability
measures, we have 
$$\left| \int_X\tilde F_n(x,1)\mu_n(x)-\int_X\tilde F(x,1)\mu_n(x)\right|
\leq \lVert \tilde F_n-\tilde F\rVert_\infty\longrightarrow 0.$$
Since $\mu_n\to \m$ weakly, we also have 
$$\left|\int_X\tilde F(x,1)\mu_n(x)-\int_X\tilde F(x,1)\mu(x)\right|
\longrightarrow 0.$$
Therefore ${\cal S}(\mu_n,\phi^n_t)([f])=\int_X\tilde F_n(x,1)\mu_n(x)\to 
{\cal S}(\mu,\phi_t)([f])=\int_X\tilde F(x,1)\mu_0(x)$.
\end{Proof}
\begin{Def}{\bf[Schwartzman unique ergodicity]} \rm We say that a flow $\phi_t$
is Schwartzman uniquely ergodic if ${\cal S}(\phi_t)$ is reduced to one point.
\end{Def}

\begin{Teo} The set ${\mathfrak S}(N)$ of Schwartzman uniquely ergodic flows 
is a $G_\delta$-set in ${\mathfrak F}(X)$.
\end{Teo}
\begin{Proof}. Fix some norm on $\rH^1(X;\R)$. We shall measure diameters of 
subsets of 
$\rH^1(X;\R)$. with respect to that norm. Fix $\epsilon>0$. Call 
${\cal U}_\epsilon$ the set of flows $\phi_t$ such that the diameter of 
${\cal S}(\phi_t)\subset \rH_1(X;\R)$ is 
$<\epsilon$. If $\phi^0_t\in {\cal U}_\epsilon$, we can find $U$ an open subset
of $\rH_1(X;\R)$ of diameter $<\epsilon$ and containing ${\cal S}(\phi^0_t)$. 
By the lemma above the set $\{\phi_t\in {\mathfrak F}(X)\mid {\cal S}(\phi_t)
\subset U\}$ is open in ${\mathfrak F}(X)$ contains $\phi^0_t$ and is 
contained in ${\cal U}_\epsilon$.
The set of Schwartzman uniquely ergodic flows is $\cap_{n\geq 1}
{\cal U}_{1/n}$.\end{Proof}

\noindent{\bf Examples.} By the computation done above, linear flows on the 
torus $\T^n$ are Schwartzman uniquely ergodic. Of course, all uniquely ergodic 
flows (\ie flows having 
exactly one invariant probability measure) are also Schwartzman uniquely 
ergodic. Moreover, by Corollary \ref{SchCon}, any flow topologically conjugate
to a Schwartzman uniquely ergodic flow is itself Schwartzman uniquely ergodic.
\\
Moreover, other examples can be obtained by the following result.

\begin{Prop}
Let $\f_t:X\longrightarrow X$ be a continuous flow on the compact path 
connected space $X$. Suppose that there exist 
$t_i\uparrow +\infty$ such that 
$\f_{t_i}{\longrightarrow} \f$ in $C(X,X)$ (with the $C^0$-topology). Then, 
$\f_t$ is Schwartzman uniquely ergodic. 
In particular, periodic flows and (uniformly) recurrent flows are Schwartzman 
uniquely 
ergodic (in both cases $\f={\rm Id}$).
\end{Prop}
\begin{Proof} Fix a continuous map $f:X\to\T$. Consider the function 
$F:X\times [0,+\infty)\to \T, (x,t)\mapsto f(\phi_t(x))-f(x)$. We have 
$F(x,0)=0$, for every $x\in X$. 
Call $\bar F:X\times [0,+\infty)\to \R$ the (unique) continuous lift of 
$F$ such that
$\bar F(x,0)=0$, for every $x\in X$. The definition of the Schwartzman 
asymptotic cycle gives
$${\cal S}(\mu)([f])=\int_X\bar F(x,1)\,d\mu(x),$$
for every probability measure invariant under $\phi_t$.
We claim that we have
$$\forall t,t'\geq 0,\forall x\in X, \quad \bar F(x,t+t')=\bar F(\phi_t(x),t')
+F(x,t).$$
In fact, if we fix $t$ and we consider each side of the equality above as a 
(continuous) function of $(x,t')$ with values in $\R$, we see that the two 
sides are equal for $t'=0$, and that they both lift the function
$$(x,t')\mapsto f(\phi_{t+t'}(x))-f(x)=f(\phi_t'(\phi_t(x))-f(\phi_t(x))+
f(\phi_t(x))-f(x)$$
with values in $\T$.
By induction, it follows easily that 
$$\forall k\in\N, \quad \bar F(x, k)=\sum_{j=0}^{k-1}\bar F(\phi_j(x),1).$$
Therefore, if $t\geq 0$ and $[t]$ is its integer part, we also obtain
\begin{equation*}
\bar F(x,t)=\bar F(\phi_{[t]}(x),t-[t])+\sum_{j=0}^{[t]-1}\bar F(\phi_j(x),1).
\tag{$*$}
\end{equation*}
It follows that
\begin{equation*}
\forall t\geq 0,\forall x\in X, \quad \lvert \bar F(x,t)\rvert\leq 
([t]+1)\lVert \bar F|X\times [0,1]\rVert_\infty.\tag{$**$}
\end{equation*}
By compactness $\lVert \bar F|X\times [0,1]\rVert_\infty$ is finite.
If we integrate equality ($*$) with respect to a probability measure $\mu$ on 
$X$ invariant under the flow $\phi_t$,  we obtain
$$\int_X\bar F(x,t)\,d\mu(x)=\int_X\bar F(x,t-[t])\,d\mu(x)+[t]\int_X\bar 
F(x,1)\,d\mu(x).$$
Therefore we have
\begin{equation*}
{\cal S}(\mu)([f])=\lim_{t\to+\infty}\int_X\frac{\bar F(x,t)}t\,d\mu(x).
\tag{$***$}
\end{equation*}
If $\gamma: [a,b]\to\T$ is a continuous path, with $a\leq b$, denote
${\cal V}(\gamma |[a,b]):=\bar\gamma(b)-\bar\gamma(a) $, where
$\bar \gamma: [a,b]\to\R$  is a continuous lift of $\gamma$ 
(this quantity does not depend on the chosen lift).\\
Suppose now that we set $\gamma_x(s)=\phi_s(x)$; we have 
$\bar F(x,t)={\cal V}(f\gamma_x|[0,t])$.
Fix now some point $x_0\in X$, and consider $t_i\to +\infty$ such that 
$\phi_{t_i}\to \phi$
in the $C^0$ topology. Since $\bar F(x_0,t)/t$ is bounded in absolute value by
$2 \lVert \bar F|X\times [0,1]\rVert_\infty$, for $t\geq 1$, extracting a 
subsequence if necessary, we can assume that $\bar F(x_0,t_i)/t_i\to c\in \R$.
If $x\in X$, we can find a continuous path $\gamma:[0,1]\to M$ with 
$\gamma(0)=x_0$
and $\gamma(1)=x$. The map $\Gamma:[0,1]\times [0,t]\to \T,(s,s')\mapsto 
\phi_{s'}(\gamma(s))$ is continuous, therefore we can lift it to a continuous 
function with values in $\R$, and this implies the equality
$${\cal V}(\Gamma|[0,1]\times\{0\})+{\cal V}(\Gamma|\{1\}\times[0,t])
-{\cal V}(\Gamma|[0,1]\times\{1\})-{\cal V}(\Gamma|\{0\}\times[0,t])=0.$$
This can be rewritten as
$${\cal V}(f\gamma_x|[0,t])-{\cal V}(f\gamma_{x_0}|[0,t])={\cal V}
(f\phi_t\gamma)-{\cal V}(f\gamma),$$
which translates to
$$\bar F(x,t)-\bar F(x_0,t)={\cal V}(f\phi_t\gamma)-{\cal V}(f\gamma).$$
Since $\phi_{t_i}\to \phi$ uniformly, by continuity of ${\cal V}$, the left 
hand-side remains bounded as $t=t_i\to +\infty$. It follows that 
$(\bar F(x,t_i)-\bar F(x_0,t_i))/t_i\to 0$.
Hence for every $x\in X$, we also have that $\bar F(x,t_i)/t_i$ tends to the 
same limit $c$ as $\bar F(x_0,t_i)/t_i$. Since $\bar F(x,t)/t$ is uniformly 
bounded for $t\geq 1$,
by $(**)$, by Lebesgue's dominated convergence we obtain from ($***$) that
${\cal S}(\mu)([f])=c$, where $c$ is independent of the invariant measure 
$\mu$. This is of course true for any $f:X\to \T$. Therefore ${\cal S}(\mu)$ 
does not depend on the invariant measure $\mu$.
 \end{Proof}
An interesting property of Schwartzman uniquely ergodic flows (which also 
shows that they have some kind of rigidity) is the following 
proposition, that follows immediately from the definition of Schwartzman unique
ergodicity and what we remarked, in the examples above, about the asymptotic 
cycles of fixed and periodic points (see also \cite{Schwartzman}).
 
\begin{Prop}
Suppose that $\phi_t $ is a Schwartzman uniquely ergodic flow on $X$. If there 
exists either a fixed point or a closed orbit homologous to zero, then all 
closed orbits are homologous to zero. In the remaining case, if $C_1$ and 
$C_2$  are closed orbits with periods $\t_1$ and $\t_2$, then $\frac{C_1}
{\t_1}$ and $\frac{C_2}{\t_2}$ are homologous. Since $[C_1]$ and $[C_2]$ are 
in $\rH_1(X;\Z)$, it follows in this case that the ratio of the periods of any 
two closed orbits must be rational. Consequently, for any continuous family of 
periodic orbits of $\f_t$, all orbits have the same period.
\end{Prop}  

\begin{Def}{\bf[Schwartzman strict ergodicity]} \rm We say that a flow 
$\phi_t$ is Schwartzman strictly ergodic if it is  Schwartzman uniquely 
ergodic and it has an invariant measure $\mu$ of full support (\ie $\mu(U)>0$ 
for every non-empty open subset $U$ of $X$).
\end{Def}

\noindent{\bf Examples.}
Linear flows on the torus $\T^n$ are Schwartzman strictly ergodic (they 
preserve Lebesgue measure). Of course, all strictly ergodic flows (\ie flows 
having exactly one invariant probability measure, and the support of this 
measure is full) are also Schwartzman strictly ergodic. A minimal flow which 
is  Schwartzman uniquely ergodic is in fact Schwartzman strictly ergodic 
(because all invariant measures have full support). Moreover, any flow 
topologically conjugate 
to a Schwartzman strictly ergodic flow is also  Schwartzman strictly ergodic.

%%%%%%%%%%%%%%%%%%%REFERENCES %%%%%%%%%%%%%%%%%%%%%%%%%%%%%%%%%%%%%%%%%%%%%%%%

\def\cprime{$'$}

\vspace{1.truecm}

{\sc Albert Fathi}\\
{\it\'Ecole normale sup\'erieure de Lyon,\\
Unit\'e de Math\'ematiques Pures et Appliqu\'ees, UMR CNRS 5669,\\
46, all\'ee d'Italie,\\
69364 Lyon Cedex 07, France}\\

{\sc Alessandro Giuliani}\\
{\it Dipartimento di Matematica,\\
Universit\`a degli Studi di Roma Tre,\\
L.go S. Leonardo Murialdo 1,\\
00146 Roma, Italy}\\

{\sc Alfonso Sorrentino}\\
{\it CEREMADE, UMR CNRS 7534,\\
Universit\'e Paris-Dauphine,\\
Pl. du Maréchal De Lattre De Tassigny,\\
75775 Paris Cedex 16, France.}

\end{document}